\documentclass[a4paper]{article}

\usepackage{enumerate,amsmath,amssymb,amsfonts,amsthm}
\usepackage{url,hyperref}
\usepackage[english]{babel}
\usepackage[all]{xy}
\usepackage{times}
\usepackage[flushmargin]{footmisc}

\newtheorem{lemma}{Lemma}
\newtheorem{proposition}[lemma]{Proposition}
\newtheorem{corollary}[lemma]{Corollary}
\newtheorem{definition}[lemma]{Definition}
\newtheorem{question}[lemma]{Question}
\newtheorem{theorem}[lemma]{Theorem}
\newtheorem{remark}[lemma]{Remark}

\makeatletter
\def\blfootnote{\xdef\@thefnmark{}\@footnotetext}
\makeatother

\newcommand{\bs}{\backslash}
\DeclareMathOperator{\Ext}{Ext}

\newcommand{\Div}{\ensuremath{\text{Div}}}

\newcommand{\Pic}{\ensuremath{\text{Pic}}}
\newcommand{\barx}{\ensuremath{\overline{X}}}
\newcommand{\bary}{\ensuremath{\overline{Y}}}

\newcommand{\Gal}{\ensuremath{\text{Gal}}}
\newcommand{\Hom}{\ensuremath{\text{Hom}}}
\newcommand{\id}{\ensuremath{\text{id}}}

\newcommand{\ob}{\ensuremath{\text{ob}}}

\newcommand{\barr}{\ensuremath{\overline{\Rcal}}}
\newcommand{\Ind}{\ensuremath{\text{Ind}}}

\newcommand{\kbar}{\ensuremath{\bar{k}}}
\newcommand{\bark}{\ensuremath{\bar{k}}}
\newcommand{\bPic}{\ensuremath{\textbf{Pic}}}
\newcommand{\cPic}{\ensuremath{\mathcal{P}\text{ic}}}
\newcommand{\coker}{\ensuremath{\text{coker}}}
\newcommand{\kar}{\ensuremath{\text{char}}}
\newcommand{\barxp}{\ensuremath{\overline{X'}}}

\newcommand{\Rcal}{\ensuremath{\mathcal{R }}}

\newcommand{\Cb}{\ensuremath{\mathbb{C }}}

\title{The elementary obstruction and the Weil restriction}
\author{Tim Wouters}

\begin{document}
\maketitle

\blfootnote{
K.U.Leuven - Departement Wiskunde, Celestijnenlaan 200B bus 2400, B-3001 Heverlee, Belgium, tim@wouters.in
\\ \ 
\\ \ \\\textit{Dedicated to the memory of Joost van Hamel, for all trouble he took to support 
this article with discussions, comments and suggestions at all times, 
also in moments when it is not obvious that one keeps on focusing on mathematics.}}

\begin{abstract}
In this text we investigate the good behaviour of the elementary
obstruction, introduced by Colliot-Th\'el\`ene and Sansuc
\cite{collsans}.  This is an obstruction to the existence of a
rational points on certain algebraic varieties. Assuming some
conditions on the Picard group, we prove that the elementary
obstruction behaves well under the Weil restriction of a variety.

\end{abstract}

\vspace{1cm}

\section{Introduction}

For a field $k$ and a variety $X$ over $k$ (i.e. a separated $k$-scheme of finite type), 
questions concerning $k$-rational points of $X$ have been studied since ages.
Different aspects arise in this area of research.
In this paper we will focus on a certain obstruction to the
existence of a rational point, namely the elementary obstruction, introduced by
Colliot-Th\'el\`ene and Sansuc \cite[Sec. 2.2]{collsans}.

Let $\kbar$ be a separable closure of $k$ and $\Gamma=\Gal(\kbar/k)$.
If $X$ is a smooth, geometrically integral variety over $k$, 
the elementary obstruction $\ob(X)$ of $X$ is defined as the class of the exact
sequence of left $\Gamma$-modules
\begin{equation*}
\text{OB}(X):= \quad 1 \to \bark^\times \to \bark(X)^\times \to \bark (X)^\times/\bark^\times \to 1
\end{equation*}
in $\Ext_\Gamma^1(\bark (X)^\times/\bark ^\times, \bark^\times)$.  Note that we use 
the common notation $\kbar(X)$ for the function field of $\barx=X \times_k \kbar$.
Analogously, we will denote $\bark [X]$ to be ring of regular functions on $\barx$.
If $X$ contains a $k$-rational point, then $\ob(X)=0$ \cite[Prop. 2.2.2]{collsans}.
Furthermore, if $\bark [X]^\times = \bark^\times$, the class of
\[ E(X):= \quad 1 \to \kbar^\times \to \kbar (X)^\times \to \Div(\barx) \to \Pic(\barx) \to 1 \]
in $\Ext_\Gamma^2 (\Pic(\barx),\kbar^\times)$ is denoted by $e(X)$.  Colliot-Th\'el\`ene and
Sansuc showed that the morphism
\[ \delta : \Ext_\Gamma^1 (\kbar(X)^\times /\kbar^\times, \kbar^\times) \to
\Ext_\Gamma^2 (\Pic (\barx),\kbar^\times), \]
which arises in the long exact sequence induced by
\[ 1 \to \kbar (X) ^\times/\kbar^\times \to \Div(\barx) \to \Pic(\barx) \to 1, \]
is injective and that $\delta(\ob (X))=e(X)$ \cite[Prop. 2.2.4]{collsans}.  This is a consequence of
Shapiro's Lemma and Hilbert 90.  Therefore it is also justified to say
$e(X)$ is the elementary obstruction of $X$.  In this paper we will mainly
use this definition for the elementary obstruction.

Several authors have been wondering whether the elementary obstruction
behaves well under classical geometric constructions.  A first observation is that the elementary obstruction is a birational invariant,
since birationally equivalent varieties have isomorphic function fields.
Wittenberg proved being zero behaves well under rational maps
 \cite[Lem. 3.1.2]{witalbelem}.
Borovoi, Colliot-Th\'el\`ene and Skorobogatov wondered whether being zero behaves
well under base extension (i.e. whether $\ob(X)=0$
implies $\ob(X \times_k K)=0$ for $K\supset k$ and $X$ a smooth, geometrically integral
variety over $k$)
\cite[Sec. 2]{borcolsko}. They gave several (partial) positive
answers to this question. Wittenberg gave a positive answer to this question for
arbitrary (smooth, proper, geometrically integral) $X$ when $K$ is a $p$-adic or real closed field \cite[Cor.
3.2.3]{witalbelem} or when $k$ is a number field and the
Tate-Shafarevich group of the Picard variety of $X$ is finite
\cite[Cor. 3.3.2]{witalbelem}.  But recently he gave a negative answer to this question 
by producing a counterexample over $\Cb((t))$ (not published yet).

In this paper we focus on the question whether being zero
behaves well under the Weil restriction of varieties.
To describe the problem more explicitely, we first recall the definition of
the Weil restriction.

\begin{definition}
Let $k$ be a field and $k'$ a finite field extension of $k$.  Let $X$ be
a variety defined over $k'$.  We say a variety $\mathcal{R}_{k'/k} X$ over $k$
is the Weil restriction of $X$ if there is a $k'$-morphism
$\varphi: \mathcal{R}_{k'/k} X \times_k k' \to X$ such that
for any  $k$-variety $Y$ and  $k'$-morphism
$f:Y\times_k k' \to X$, a unique
$k$-morphism $g:Y\to \mathcal{R}_{k'/k} X$ exists such that
$\varphi \circ g' = f$, where $g': Y\times_k k' \to \mathcal{R}_{k'/k} X \times_k k'$
 is the $k'$-morphism induced by $g$.
If the Weil restriction exists, it is unique up to $k$-isomorphism.
\end{definition}

\noindent There is a well known proposition that guarantees the existence of the Weil restriction
under the assumption of some conditions.

\begin{proposition} \label{prop:existweil}
Let $k$ be a field, $\kbar$ a separable closure and $k'$ a finite subextension of $k$ in $\kbar$.
Denote $\Gamma=\Gal(\kbar/k)$, $H=\Gal(\kbar/k')$ and let $X$ be a quasiprojective variety over $k'$. The Weil restriction
$\mathcal{R}_{k'/k} X$ of $X$ exists, and
\[ \mathcal{R}_{k'/k} X \times_{k} \kbar = \prod_{[\sigma] \in H \bs \Gamma } \sigma X \]
where $\sigma X$ is the $\kbar$-variety obtained by base extension of $X\times_k \kbar$ by $\sigma:\kbar \to \kbar$ and $H \bs \Gamma $ are the right cosets of $H$ in $\Gamma$.
The $k'$-morphism $\varphi: \mathcal{R}_{k'/k} X \times_k k' \to X$ is obtained by descent theory from its base
extension
$\overline{\varphi}: \overline{ \mathcal{R}_{k'/k} X } \to \overline{X}$, which is the projection
onto the factor $(\id ) X$.
\end{proposition}

\noindent For the proof we refer to \cite[Prop. 16.26]{milnealggeom}.  Remark that if
$[\sigma]=[\tau]\in \Gamma/H$, the universal property of fibre products
garantuees
$\sigma X$ and $\tau X$ to be isomorphic as $\kbar$-varieties.
The universal property of the Weil restriction implies also a 1-1 correspondence
between $\mathcal{R}_{k'/k} X (k)$ and $X(k')$, since rational points
are equivalent with sections of the structure morphism.  As there is a natural connection between rational
points of a variety and its Weil restriction, it is natural to
ask the following question.

\begin{question}
Let $k$ be a field and $k'$ a finite field extension.  Suppose $X$ is a smooth, geometrically integral variety
over $k'$ such that the Weil restriction $\mathcal{R}_{k'/k}X$ exists.  Does
$e(X)=0$ implies $e(\mathcal{R}_{k'/k}X)=0$ and vice versa?
\end{question}

We will answer this question partially positively (Proposition \ref{prop:weilelobstr} and Theorem \ref{prop:elobstrweilcond}),
but first we will give a result on the elementary obstruction of a product variety (Theorem \ref{prop:prodelobstr}),
as Proposition \ref{prop:existweil} tells us that the Weil restriction is
closely related to product varieties.

\section{Product varieties}

In this section let $k$ be a field, $\bar{k}$ a separable closure and
$\Gamma=\Gal( \bar{k}/k)$.
Let $X$ and $Y$ be two smooth geometrically integral varieties over $k$, then
the following theorem is a merely homological result.

\begin{theorem} \label{prop:prodelobstr}
The multiplication $\Gamma$-morphism $\pi:\bar{k}(X)^\times/\bar{k}^\times \oplus \bar{k} (Y)^\times/\bar{k}^\times \to \bar{k}(X\times_k Y)^\times/\bar{k}^\times$
induces a morphism by pullback
\begin{multline*}
 {\pi^\ast }': \Ext^1_\Gamma (\bar{k}(X\times_k Y)^\times/\bar{k}^\times,\bar{k}^\times ) \to \\
\Ext^1_\Gamma (\bar{k}(X)^\times/\bar{k}^\times,\bar{k}^\times ) \oplus
\Ext^1_\Gamma (\bar{k}(Y)^\times/\bar{k}^\times,\bar{k}^\times ) 
\end{multline*}
such that ${\pi^\ast}' (\ob (X \times_k Y))=(\ob (X),\ob (Y))$.
If $\bar{k}[X]^\times = \bar{k}^\times =  \bar{k}[Y]^\times $,
then the $\Gamma$-morphism $\psi:\Pic(\overline{X}) \oplus \Pic(\overline{Y}) \to \Pic(\overline{X}\times_{\bar{k}} \overline{Y})$
defined by pullback of line bundles, induces a morphism
\[ {\psi^\ast}':
\Ext^2_\Gamma (\Pic(\overline{X}\times_{\bar{k}} \overline{Y}),\bar{k}^\times ) \to
\Ext^2_\Gamma (\Pic(\overline{X}),\bar{k}^\times ) \oplus
\Ext^2_\Gamma (\Pic(\overline{Y}),\bar{k}^\times ) \]
such that ${\psi^\ast}'(e (X \times_k Y))=(e (X),e (Y))$.
Even more ${\pi^\ast}'$ and ${\psi^\ast}'$ commute with the natural inclusions
in the following commutative diagram:
\[ \xy
\xymatrix"*"{%
 \Ext_\Gamma^1 (\bar{k}(Y)^\times /\bar{k}^\times ,\bar{k}^\times) \oplus%
 \Ext_\Gamma^1  (\bar{k}(Y)^\times /\bar{k}^\times ,\bar{k}^\times)  \ar[d]^{\delta}  \\%
\Ext_\Gamma^2 (\Pic(\barx) ,\bar{k}^\times) \oplus \Ext_\Gamma^2 (\Pic(\bary) ,\bar{k}^\times).%
}
\POS-(50,-8)
\xymatrix{
\Ext_\Gamma^1 (\bar{k}(X\times_k Y)^\times /\bar{k}^\times ,\bar{k}^\times) \ar[d]^{\delta} \ar["*"]^{{\pi^\ast}'}  \\ 
\Ext_\Gamma^2 (\Pic(\barx \times_{\kbar} \bary) ,\bar{k}^\times) \ar["*"]^{{\psi^\ast}'} \\ 
}
\endxy 
\]
If $\pi$ or $\psi$ is an isomorphism, then $e(X\times_k Y)=0$ (resp. $\ob (X\times_k Y)=0$) if
and only if $e(X)=0$ and $e(Y)=0$ (resp. $\ob(X)=0$ and $\ob (Y)=0$).
\end{theorem}

\begin{remark} \label{rem:prodelobstr}
Note that if $X$ and $Y$ are smooth geometrically integral varieties satisfying
$\bar{k}[X]^\times=\bar{k}^\times=\bar{k}[Y]^\times$, then $X\times_k Y$ is also smooth
geometrically integral, and by a result of Rosenlicht \cite[Thm. 2]{rosentoralggr} it satisfies $\bar{k}[X\times_k Y]^\times =\kbar^\times$.  So speaking about $e (X \times_k Y)$ in the
second case does make sense.
\end{remark}

\begin{proof}
If we denote the canonical isomorphism
\begin{multline*} \Ext^1_\Gamma (\bar{k}(X)^\times/\bar{k}^\times \oplus \bar{k} (Y)^\times/\bar{k}^\times,\bar{k}^\times)\cong \\
\Ext^1_\Gamma (\bar{k}(X)^\times/\bar{k}^\times ,\bar{k}^\times) \oplus
\Ext^1_\Gamma ( \bar{k} (Y)^\times/\bar{k}^\times,\bar{k}^\times) \end{multline*}
by $\varphi$, then ${\pi^\ast}':=\varphi \circ {\pi^\ast}$ is the required morphism, where
\[ \pi^\ast:\Ext^1_\Gamma(\bark(X\times_k Y)^\times/\bark^\times ,\bark^\times) \to
\Ext^1_\Gamma (\bar{k}(X)^\times/\bar{k}^\times \oplus \bar{k} (Y)^\times/\bar{k}^\times,\bar{k}^\times)\]
is the pullback of 1-extensions by $\pi$.  We now prove the assertion
on the elementary obstruction.

We surely have a morphism of short exact sequences which consists of product morphisms:
\[
\xymatrix{
 1 \ar[d]  & 1 \ar[d] \\
\bar{k}^\times \oplus \bar{k}^\times \ar[d] \ar[r]^{\pi_1} 
& \bar{k}^\times \ar[d]  \\
\bar{k}(X)^\times \oplus \bar{k} (Y)^\times \ar[d] \ar[r]^{\pi_2} 
& \bar{k}(X\times_k Y)^\times \ar[d] \\
\bar{k}(X)^\times/\bar{k}^\times \oplus \bar{k} (Y)^\times/\bar{k}^\times \ar[d] \ar[r]^{\quad \ \ \pi_3=\pi} 
& \bar{k}(X\times_k Y)^\times/\bar{k}^\times \ar[d]  \\
1 & 1. }
\]Denote the left short exact sequence  by $E(X)\oplus E(Y)$.  By the notation introduced
in the introduction, the right short exact sequence is denoted by $E(X\times_k Y)$.
By the general theory of Yoneda extensions \cite[Ch. III]{maclane}, we get
\[ \varphi^{-1} (e(X),e(Y))= [\pi_1 (E(X)\oplus E(Y))] = [E(X\times_k Y) \pi_3] = \pi^\ast(e(X\times_k Y)),\]
where $\pi_1 (E(X)\oplus E(Y))$ denotes the pushforward of the
Yoneda extension $E(X)\oplus E(Y)$ by $\pi_1$ and $E(X\times_k Y) \pi_3$
denotes the pullback of the Yoneda extension $E(X\times_k Y)$ by $\pi_3$. This
proves the first part.

The second part is proved analogously, using $\Gamma$-morphisms
$\pi_4:\Div(\overline{X}) \oplus \Div(\overline{Y}) \to \Div(\overline{X}\times_{\bar{k}} \overline{Y})$
and
$\psi:\Pic(\overline{X}) \oplus \Pic(\overline{Y}) \to \Pic(\overline{X}\times_{\bar{k}} \overline{Y})$.
The commutativity assertion follows from the following morphism of short exact sequences

\[
\xymatrix{
1 \ar[d] &  1 \ar[d] \\
\bar{k}(X)^\times/\bar{k}^\times \oplus \bar{k}(Y)^\times /\bar{k}^\times \ar[d] \ar[r]^{\pi_3} 
& \bar{k}(X\times_k Y)^\times/\bar{k}^\times \ar[d] \\
\Div(\barx) \oplus \Div(\bary) \ar[d] \ar[r]^{\pi_4} 
& \Div(\barx \times_{\kbar} \bary) \ar[d]  \\
\Pic(\barx) \oplus \Pic(\bary) \ar[d] \ar[r]^{\pi_5=\psi} 
& \Pic(\barx \times_{\kbar} \bary) \ar[d] \\
1 & 1. \\
}
\]
which induces a morphism of long exact sequences, by Shapiro's lemma and Hilbert 90
containing the required diagram.

So we see that in any case $e(X)=0$ and $e(Y)=0$ (resp. $\ob(X)=0$ and $\ob(Y)=0$) if 
$e(X\times Y)=0$ (resp. $\ob(X\times Y)=0$).  If $\psi$ (resp. $\pi$) is an isomorphism,
${\psi^\ast}'$ (resp. ${\pi^\ast}'$) will be so too, so in one of these cases the inverse 
implication holds as well (recall that $e(-)=0$ if and only if $\ob(-)=0$).  
\end{proof}

\begin{remark} \label{rem:dirsum}
A known result says that if $\barx$ and $\bary$ are varieties over separable closed field $\kbar$, then as groups the morphism $\psi: \Pic(\barx) \oplus \Pic (\bary) \to \Pic(\barx\times_{\kbar} \bary)$,
defined by pull-backs, has a section.
This section restricts a line bundle
on $X\times_k Y$ to ${x_0}\times Y$ and $X\times {y_0}$ where $x_0$ and $y_0$ are
base points on $X$ and $Y$.  So as groups $\Pic(\barx)\oplus \Pic(\bary)$ is a direct
summand of $\Pic(\barx\times_{\kbar} \bary)$.  This looks interesting  to
get more information on the structure of $\Ext_\Gamma^2(\Pic(\barx\times_{\kbar} \bary),\kbar^\times)$.
In our case however, $X$ and $Y$ are
defined over a not necessarily separably closed field $k$ and 
$ \psi: \Pic (\barx)\oplus \Pic(\bary) \to
\Pic(\barx\times_{\kbar} \bary) $
 is
a $\Gamma$-morphism, but the section is not necessarily a $\Gamma$-morphism since the base points do not
have to behave well (if we do not know anything about the existence of $k$-rational points
on $X$ and $Y$).  So we can not use this result to extend the previous theorem in a direct way.  However,
we do retrieve $\psi$ is injective.
\end{remark}

Off course $\psi:\Pic(\overline{X}) \oplus \Pic(\overline{Y}) \to \Pic(\overline{X}\times_{\bar{k}} \overline{Y})$ does not need to be an isomorphism,
the product of an elliptic curve with itself delivering a counterexample
\cite[Ch. IV, Ex. 4.10]{hart}.  We can however give sufficient
conditions for $\psi$ to be an isomorphism.   This will involve the notion
of the \emph{relative Picard functor} and the \emph{Picard variety}.  If $X$ is a smooth,
geometrically integral, projective  variety over a field
$k$, we denote the relative Picard functor by $\cPic_{X/k}$, which is
representable by a group variety $\bPic (X)$, the Picard variety.
Denote by $\bPic^0 (X)$ the zerocomponent of $\bPic (X)$.  (See \cite[Ch. 8]{neron} for more information.)

\begin{proposition} \label{prop:pic00}
If $X$ is projective and $\bPic^0 (\barx)=0$, then
$\psi: \Pic(\barx) \oplus \Pic (\bary) \to \Pic(\barx\times_{\bark} \bary)$ is a
$\Gamma$-isomorphism.
\end{proposition}
\begin{proof}
By Remark \ref{rem:dirsum} we know that $\psi$ is injective, so 
it is sufficient to
prove $\coker \ \psi =0$.
By definition
\[ \cPic_{\barx/\bark} (\bary) = \Pic (\barx \times_{\bark} \bary)/\Pic(\bary) \cong \Hom_{\bark} (\bary ,\bPic (\barx)).\]
Any $f\in \Hom_{\bark} (\bary,\bPic (\barx))$ has a connected image, but
since $\bPic^0 (\barx)=0$, the connected components of $\bPic (\barx)$ are its points.  So
$\Hom_{\bark} (\bary ,\bPic (\barx))$ consists of the constant maps onto a point of $\bPic (\barx)$.
This does not depend on $Y$, so
\[ \Hom_{\bark} (\bary,\bPic (\barx))\cong \Hom_{\bark} (\bark,\bPic (\barx)) \cong \Pic(\barx). \]
Because these isomorphisms are induced by the representability
of the Picard functor,
\[ \coker \ \psi = \frac{\Pic (\barx \times_{\bark} \bary)/\Pic(\bary)}{\Pic (\barx)}
= \frac{\Pic(\barx)}{\Pic (\barx)}=0. \]
\end{proof}

\begin{proposition} \label{prop:picfin}
If $X$ is quasiprojective, $\kar (k)=0$ and $\Pic(\barx)$ is finitely generated,
then $\Pic(\barx) \oplus \Pic (\bary) \cong \Pic(\barx\times_{\bark} \bary)$.
\end{proposition}
\begin{proof}
Say $X\subset X_1$ for a projective variety $X_1$.  Since $\kar (k)=0$, there
exists a (smooth, projective) Hironaka desingularisation $X'$  of $X_1$.
As $X$ is smooth, $X$ is isomorphic to an open of $X'$.  So without
loss of generality we assume $X$ is an open part of $X'$.
The exact sequence
\[ \Div_{\barxp \bs \barx}(\barxp) \to \Pic(\barxp) \to \Pic(\barx) \to 0 \]
induces $\Pic(\barxp)$ to be finitely generated, as
$\Pic(\barx)$ and $\Div_{\barxp \bs \barx}(\barx)$ are finitely generated.  ($\Div_{\barxp \bs \barx}(\barx)$
are the divisors on $\barxp$ with support outside $\barx$.)

It suffices to prove $\Pic(\barxp \times_{\kbar} \bary)\cong \Pic(\barxp)\oplus \Pic(\bary)$
as this will induce $\Pic(\barx \times_{\kbar} \bary)\cong \Pic(\barx)\oplus \Pic(\bary)$.
Indeed, there is a commutative diagram
\[
\xymatrix{
0 \ar[r] & \Pic ( \barxp) \oplus \Pic (\bary ) \ar[r] \ar[d] &\Pic ( \barxp \times_{\bark} \bary ) \ar[d] \\
0 \ar[r] & \Pic ( \barx) \oplus \Pic (\bary ) \ar[r] \ar[d] & \Pic ( \barx \times_{\bark} \bary ) \ar[d] \\
& 0 & 0 }
\]
where the vertical arrows are the surjective restriction morphisms.  If the
injection of the first row turns out to be an isomorphism, then the injection of the
bottom row should also be surjective, hence an isomorphism.

Because $\Pic (\barxp)$ is finitely generated, we have $\bPic^0(\barxp)=0$.  
Indeed, if $\bPic^0(\barxp)\neq 0$, then $\bPic^0(\barxp)$ is an abelian variety of dimension $m>0$
 whose group of $\bark$-points is finitely generated as 
$\Pic(\barxp) = \Hom_{\kbar}(\kbar,\bPic(\barxp))$ is finitely generated. 
On the other hand the group of $\bark$-points of an abelian variety is divisible \cite[Thm. 2]{frey}.  But a divisible,
non-trivial, finitely generated group does not exist.  In this way we get a contradiction and 
so the proposition follows by Proposition \ref{prop:pic00}.
\end{proof}

\noindent Consequently we obtain the following result.

\begin{corollary}
Let $X$ and $Y$ be a smooth, geometrically integral varieties over a field
$k$ with $\bar{k}[X]^\times = \bar{k}^\times =  \bar{k}[Y]^\times $.
Let $\bark$ be a separable closure of $k$ and $\Gamma=\Gal(\bark/k)$.
If one of the following
conditions holds
\begin{enumerate}[(i)]
\item $X$ is projective and $\bPic^0 (\barx)=0$, or
\item $X$ is quasiprojective, $\kar (k)=0$ and $\Pic(\barx)$ is finitely generated,
\end{enumerate}
then
\[ {\psi^\ast}':
\Ext^2_\Gamma (\Pic(\barx \times_{\bark} \bary),\bark^\times) \to
\Ext^2_\Gamma (\Pic(\barx),\bark^\times) \oplus
\Ext^2_\Gamma (\Pic(\bary),\bark^\times)
\]
is an isomorphism such that ${\psi^\ast}' (e(X\times_k Y))=(e(X),e(Y))$.
\end{corollary}

\noindent So if one of the conditions is true, $e(X\times_k Y)=0$ if and only if $e(X)=0$ and
$e(Y)=0$.

\section{Weil restriction} \label{sec:weil}

Knowing more on the case of product varieties, we proceed to the  Weil restriction.
Throughout this section we will assume $k$ is a field, $\kbar$ is a separable closure of $k$
and $k'$ is a finite subextension  of $k$ in $\kbar$.
Denote $\Gamma=\Gal (\kbar / k)$, $H=\Gal (\kbar/k')$, and
let $X$ be a smooth, geometrically integral, quasiprojective variety over $k'$.
In this case the Weil restriction of $X$ exists by Proposition \ref{prop:existweil} and
we abbreviate it as $\Rcal$.

\begin{proposition} \label{prop:weilelobstr}
The natural $H$-morphism $\kbar(X)^\times \to \kbar(\Rcal)^\times$ induces a pullback of
1-extensions
\[ \Pi^\ast:
\Ext_\Gamma^1 (\kbar(\Rcal)^\times/\kbar^\times ,\kbar^\times) \to
\Ext_H^1 (\kbar(X)^\times / \kbar^\times,\kbar^\times),  \]
with $\Pi^\ast(\ob(\Rcal))=  \ob (X)$.  If furthermore $\bark [X]^\times=\kbar^\times$, then
the natural $H$-morphism $\Pic (\barx) \to \Pic (\barr)$ induces a pullback of 2-extensions
\[ \Phi^\ast:
\Ext_\Gamma^2 (\Pic(\barr),\kbar^\times) \to
\Ext_H^2 (\Pic(\barx),\kbar^\times) , \]
with $\Phi^\ast(e(\Rcal))=  e(X)$.
As in Theorem \ref{prop:prodelobstr} these morphisms commute
with the natural inclusions sending $\ob( - )$ to $e(- )$.
\end{proposition}

\begin{remark} \label{rem:expl}
The natural $H$-morphisms mentioned in the proposition are induced by 
Theorem \ref{prop:existweil}.  This proposition gives a 
 $k'$-morphism $\varphi:\mathcal{R}\times_k k'\to X$ retrieved by descent from
the $\kbar$-projection $\overline{\varphi}: \barr \to \barx$.  This morphism $\overline{\varphi}$
gives by pullback of principle divisors and line bundles the required $H$-morphisms.
\end{remark}

\begin{remark}
As in Remark \ref{rem:prodelobstr} it is true that
$\bark [\Rcal]^\times=\kbar^\times$ provided $\bark [X]^\times=\kbar^\times$.  So
it makes sense to speak about $e(\Rcal)$ if at first glance we only require  $\bark [X]^\times=\kbar^\times$.
\end{remark}

\begin{proof}
We give the proof of the assertion on 2-extensions.  The assertion on 1-extensions follows
in the same way, and the commutative part will follow as in Theorem \ref{prop:prodelobstr}.

Denote the $H$-morphism $\Pic (\barx) \to \Pic (\barr)$ by $\varphi'$.
This induces a pullback
\[ \varphi'^\ast: \Ext^2_H (\Pic(\barr),\kbar^\times) \to \Ext^2_H (\Pic(\barx),\kbar^\times).
\]
If we use the forgetful map
\[ \pi:\Ext^2_\Gamma (\Pic(\barr),\kbar^\times)\to \Ext^2_H (\Pic(\barr),\kbar^\times), \]
we get the required morphism $\Phi^\ast= \varphi'^\ast \circ \pi$.

To prove  $\Phi^\ast(e(\Rcal))=e(X)$, we use the following morphism of
$H$-extensions
\[
\xymatrix{
 E(X) = 
1 \ar[r]  &
\bark^\times \ar[r] \ar[d]_{\id} &
\bark (X)^\times \ar[r] \ar[d] &
\Div(\barx) \ar[r] \ar[d] &
\Pic(\barx) \ar[r] \ar[d]^{\varphi'} &
1 \\
 E(\Rcal) = 
1 \ar[r] &
\bark^\times \ar[r]  &
\bark (\Rcal)^\times \ar[r]  &
\Div(\barr) \ar[r]  &
\Pic(\barr) \ar[r]  &
1.
}
\]
As it is clear that the $H$-equivalence class of $E(\Rcal)$ equals $\pi([e(\Rcal)])$,
we get from elementary homological reasons
\[ \Phi^\ast(e(\Rcal))= \varphi'^\ast (\pi([e(\Rcal)]))= [E(X)]=e(X).\]
\end{proof}

So $e(\Rcal)=0$ implies $e(X)=0$.  We proceed
figuring out when the converse is  true.  This will hold
in the very same situation as the converse holds for product varieties.
To prove this, we will use the notion of induced group module,
with some corresponding notation. Let $G$ be a profinite group,
$H$ a subgroup of $G$ and $A$ a left $H$-module, then the induced $G$-module is
$\Ind^G_H(A):= \mathbb{Z}[G] \otimes_{\mathbb{Z}[H]} A$, where $\mathbb{Z}[G]$ is considered
as a right $\mathbb{Z}[H]$-module. 
This is a left $G$-module, the $G$-action is defined by
$\phantom{}{\gamma'}(\gamma \otimes a)= \gamma'\gamma \otimes a$ for any $a\in A$ and $\gamma,\gamma'\in G$.
If $A$ and $B$ are left $H$-modules and $f:A\to B$ is an $H$-morphism, then
we get an induced $G$-morphism
\[
\Ind^G_H (f):\Ind^G_H (A) \to \Ind^G_H (B) \quad \text{ defined by } \quad \gamma \otimes a \mapsto \gamma \otimes f(a),
\]
for any $a\in A$ and $\gamma \in G$.  If $B$ is also a left $G$-module,
we write $\Ind^G_H (f)'$ for the $G$-morphism $\pi\circ \Ind^G_H (f)$ with
\[\pi:  \Ind^G_H (B) \to B \quad \text{ defined by } \quad \gamma \otimes b \mapsto \phantom{}^\gamma b.
\]
If $E$ is an exact sequence
\[
\xymatrix{
A_{1} \ar[r]^{f_{1}} &
A_{2} \ar[r]^{f_{2}} &
A_{3},
}
\]
then we get an induced exact sequence $\Ind^G_H (E)$
\[
\xymatrix{
\Ind^G_H (A_{1}) \ar[r]^{\tilde{f}_{1}} &
\Ind^G_H (A_{2}) \ar[r]^{\tilde{f}_2} &
\Ind^G_H (A_{3}),
}
\]
where we have denoted $\tilde{f}_i:=\Ind^G_H (f_i)$ for sake of simplicity.

\begin{theorem} \label{prop:elobstrweilcond}
If $\kbar [X]^\times=\kbar^\times$ and
if one of the two following conditions is true
\begin{enumerate}[(i)]
\item \label{cond1} $X$ is projective and $\bPic^0 (\barx)=0$, or
\item \label{cond2} $X$ is quasiprojective, $\kar (k)=0$ and $\Pic(\barx)$ is finitely generated,
\end{enumerate}
then $\Phi^\ast$ of Proposition \ref{prop:weilelobstr} is an isomorphism.
\end{theorem}

\begin{proof}

We will prove this result by giving another description of $\Phi^\ast$.

If $\varphi'$ is the $H$-morphism $\Pic(\barx)\to \Pic(\barr)$ as defined
in the proof of Proposition \ref{prop:weilelobstr}, 
the induced $\Gamma$-morphism $\Ind^\Gamma_H(\varphi')':\Ind^\Gamma_H(\Pic(\barx))\to \Pic(\barr)$
gives a pullback of 2-extensions,
\begin{equation*}
 \Ind^\Gamma_H(\varphi')'^\ast: \Ext^2_\Gamma (\Pic (\barr), \kbar^\times ) \to
\Ext^2_\Gamma (\Ind^\Gamma_H \Pic (\barx), \kbar^\times ).
\end{equation*}
Furthermore say $\pi'$ is the forgetful map
\[ \pi': \Ext^2_\Gamma (\Ind^\Gamma_H(\Pic(\barx)),\bark^\times)) \to
\Ext^2_H (\Ind^\Gamma_H(\Pic(\barx)),\bark^\times)) \]
and
\[ i^\ast: \Ext^2_H (\Ind^\Gamma_H(\Pic(\barx)),\bark^\times)) \to
\Ext^2_H (\Pic(\barx),\bark^\times)) \]
the pullback by $i:\Pic(\barx)\to \Ind^\Gamma_H(\Pic(\barx)): \mathcal{L} \mapsto \id \otimes \mathcal{L}$.
So we end up in the following situation:
\[ \xy
\xymatrix"b"{
\Ext^2_\Gamma(\Pic (\barr),\bark^\times)  \ar[r]^\pi &
\Ext^2_H(\Pic (\barr),\bark^\times) \ar[r]^{\varphi'^\ast} &
\Ext^2_H (\Pic (\barx),\bark^\times)
}
\POS-(-15,15)
\xymatrix"a"{
\Ext^2_\Gamma(\Ind^\Gamma_H (\Pic (\barx)),\bark^\times ) \ar[r]^{\pi'} \ar@{<-}["b"]^(.6){\Ind^\Gamma_H(\varphi')'^\ast}
&
\Ext^2_H(\Ind^\Gamma_H (\Pic (\barx)),\bark^\times) \ar["b"r]_(.6){i^\ast}.
}
 \endxy
\]
We will prove $\Phi^\ast=\varphi'^\ast \circ \pi$ is an isomorphism by proving that $i^\ast
\circ \pi' \circ \Ind^\Gamma_H(\varphi')'^\ast$ is an isomorphism
and that the diagram above commutes.  The latter follows directly
from elementary homological reasons.

To prove the former, first observe that $i^\ast \circ \pi'$ is an
isomorphism by Shapiro's Lemma as it has an inverse $\Ind^\Gamma_H (\id)'_\ast \circ
\Ind^\Gamma_H$ with
\begin{multline*} \Ind^\Gamma_H : \Ext^2_H(\Pic(\barx),\bark^\times) \to
\Ext^2_\Gamma( \Ind^\Gamma_H (\Pic(\barx)),\Ind^\Gamma_H(\bark^\times)): \\ [E] \mapsto
[\Ind^\Gamma_H(E)] \end{multline*}
and $\Ind^\Gamma_H (\id)'_\ast$ the pushforward
\[
\Ext^2_\Gamma( \Ind^\Gamma_H
(\Pic(\barx)),\Ind^\Gamma_H(\bark^\times)) \to \Ext^2_\Gamma(
\Ind^\Gamma_H(\Pic(\barx)),\bark^\times) \] by $\Ind^\Gamma_H
(\id)':\Ind^\Gamma_H(\bark^\times)\to \bark^\times$.  This is indeed an inverse by elementary homological reasons.

So it remains to prove $\Ind^\Gamma_H(\varphi')'^\ast$ is an isomorphism.
Therefore we first
choose a set of representatives $\{ \sigma_1, \ldots ,\sigma_n \}$ of the classes of $H\bs \Gamma$
with $\sigma_1=\id$.

If Condition (\ref{cond1}) or (\ref{cond2}) is true, then pullback along all components
\[ \psi: \bigoplus_{i=1}^n \Pic(\sigma_i X) \to \Pic(\barr) \]
is an isomorphism of $H$-modules by Proposition \ref{prop:pic00} and \ref{prop:picfin}.
We will prove there is a 1-1 correspondence $\tau:\Ind^\Gamma_H(\Pic(\barx))\to
\bigoplus_{i=1}^n \Pic(\sigma_i X)$ and that
$\psi \circ \tau=\Ind^\Gamma_H(\varphi')'$, hence inducing $\Ind^\Gamma_H(\varphi')'$ is an
isomorphism.

First remark that for all $i=1,\ldots,n$, base extension by $\sigma_i$ induces
a bijection $B_i: \Pic (\barx)\to \Pic (\sigma_i X)$ which does not need to
be a $H$-morphism as $H$ does not necessarily commute with $\sigma_i$.  There are also
$H$-morphisms $\psi_i: \Pic(\sigma_i X) \to \Pic(\barr)$ induced by projection
on the $i$-th factor, so $\psi=\sum_{i=1}^n \psi_i$ and $\psi_1=\varphi'$.  It is easy to
see the $B_i$ and $\psi_i$ relate in the following way
$\phantom{}^{\sigma_i^{-1}} \psi_i(B_i(\mathcal{L}))=\psi_1(\mathcal{L})$
for any $\mathcal{L}\in \Pic(\barx)$.

To define $\tau$, it satisfies defining $\tau(\gamma \otimes \mathcal{L})$ for any
$\mathcal{L}\in \Pic(X)$ and $\gamma \in \Gamma$.
Suppose $\gamma=\sigma_i h$ for $h\in H$ and $1\leq i \leq n$,
then we set $\tau(\gamma \otimes \mathcal{L})$ with 0 as $[\sigma_j]$-components
for $j\neq i$ and  $ B_i(\phantom{}^h \mathcal{L})$ as $[\sigma_i]$-component.
This is well defined and as all the $B_i$ are bijections, $\tau$ is indeed a 1-1 correspondence.
Even more
\[\psi \circ \tau (\gamma \otimes \mathcal{L})
=\psi_i( B_i(\phantom{}^h \mathcal{L}))
=\phantom{}^{\sigma_i} \psi_1(\phantom{}^h \mathcal{L})
=\phantom{}^{\gamma} \psi_1(\mathcal{L})= \Ind^\Gamma_H(\varphi')'(\gamma \otimes \mathcal{L}). \]
\end{proof}

\noindent So if one of the two conditions holds, $e(X)=0$ if and only if $e(\Rcal)=0$.

\section*{Acknowledgement}

The author is supported by Fonds voor Wetenschappelijk Onderzoek Vlaanderen (G.0318.06).


\end{document}